\documentclass{amsart}
\usepackage[utf8]{inputenc}

\usepackage{latexsym,amssymb,amsmath}
\usepackage{epsfig}
\input xy
\xyoption{all}
\usepackage[english]{babel}

\usepackage{amsmath}
\newtheorem{theorem}{Theorem}[section]
\newtheorem{lemma}[theorem]{Lemma}
\newtheorem{proposition}[theorem]{Proposition}

\newtheorem{definition}[theorem]{Definition}

\def\PP{\mathbb{P}}
\def\CC{\mathbb{C}}
\def\ZZ{\mathbb{Z}}\def\QQ{\mathbb{Q}}\def\SS{\mathbb{S}}
\def\LL{\mathbb{L}}\def\FF{\mathbb{F}}

\def\cO{{\mathcal O}}\def\cN{{\mathcal N}}

\def\ra{\rightarrow}\def\lra{\longrightarrow}
\def\fg{{\mathfrak g}}\def\fp{{\mathfrak p}}
\def\fsl{{\mathfrak{sl}}}

\title{Topics on the geometry of homogeneous spaces}
\author{Laurent Manivel}

\address{Institut de Math\'ematiques de Toulouse ; UMR 5219 \\ %
UPS CNRS, F-31062 Toulouse Cedex 9, France}
\email{manivel@math.cnrs.fr}

\begin{document}

\begin{abstract}
This is a survey paper about a selection of results in complex algebraic geometry that 
appeared in the recent and less recent litterature, and in which rational homogeneous 
spaces play a prominent r\^ole. This selection is largely arbitrary and mainly reflects
the interests of the author. 
\end{abstract}

\maketitle

Rational homogeneous varieties are very special projective varieties, which appear in a 
variety of circumstances as exhibiting extremal behavior. In the quite recent years, a 
series of very interesting examples of pairs (sometimes called Fourier-Mukai partners) 
of derived equivalent, but not isomorphic, 
and even non birationally equivalent manifolds have been discovered by several authors,
starting from the special geometry of certain homogeneous spaces.
We will not discuss derived categories and will
not describe these derived equivalences: this would require more 
sophisticated tools and much ampler discussions. Neither will we say much about Homological
Projective Duality, which can be considered as the unifying thread of all these apparently 
disparate examples. Our much more modest goal will be to describe their geometry, starting 
from the ambient homogeneous spaces. 

In order to do so, we will have to explain how one can approach homogeneous spaces just
playing with Dynkin diagram, without knowing much about Lie theory. In particular we will 
explain how to describe the VMRT (variety of minimal rational tangents) of a generalized 
Grassmannian. This will show us how to compute the index of these varieties, remind
us of their importance in the classification problem  of Fano manifolds, in rigidity questions, 
and also, will explain their close relationships with prehomogeneous vector spaces. 

We will then consider vector bundles on homogeneous spaces, and 
use them to construct interesting birational transformations, including important types
of flops: the Atiyah and Mukai flops, their stratified versions, also the Abuaf-Segal
and Abuaf-Ueda flops; all these beautiful transformations are easily described in 
terms of homogenenous spaces. And introducing sections of the bundles involved, we 
will quickly arrive at several nice examples of Fourier-Mukai partners.

We will also explain how the problem of finding crepant resolutions of orbit closures 
in prehomogeneous spaces is related to the construction of certain manifolds
with trivial canonical class. This gives a unified perspective over  classical constructions
by Reid, Beauville-Donagi and Debarre-Voisin of abelian and hyperK\"ahler varieties,
naturally embedded into homogeneous spaces. The paper will close on a recent 
construction, made in a similar spirit, of a generalized Kummer fourfold from an 
alternating three-form in nine variables.

\smallskip\noindent {\bf Acknowledgements}. These notes were written for the Lecture Series in Algebraic Geometry, organized in August and September 2019 at the Morningside Center of Mathematics
in Beijing. We warmly thank Professors Baohua Fu, Yujiro Kawamata and Shigeru Mukai for the 
perfect organization. 

\section{Basics}

\subsection{Rational homogeneous spaces}
A classical theorem of Borel and Remmert \cite{akhiezer} 
asserts that a projective complex manifold 
which admits a transitive action of its automorphism group is a direct product of an
abelian variety by a rational homogeneous space. The latter can be described as a 
quotient $G/P$, where $G$ is a semi-simple algebraic group and $P$ a {\it parabolic
subgroup}. Moreover it can always be decomposed into a product 
$$G/P \simeq G_1/P_1 \times\cdots \times G_\ell /P_\ell$$
of rational homogeneous spaces of simple algebraic groups $G_1, \ldots , G_\ell$. 

So we will suppose in the sequel, unless otherwise stated, that $G$ is simple. 
Moreover the list of homogeneous spaces under $G$ will only depend on its Lie algebra $\fg$,
which is clasically encoded in a Dynkin diagram $\Delta$. In fact the $G$-equivalence classes 
of parabolic subgroups are in bijective correspondence with the finite subsets of nodes of
$\Delta$. As a conclusion, a projective rational homogeneous space with simple automorphism
group is determined by a marked Dynkin diagram.

The two extremes cases correspond to the {\it complete flag manifolds} (all nodes marked),
and the {\it generalized Grassmannians} (only one node marked). Generalized Grassmannians 
are equivariantly embedded inside the projectivizations of the fundamental representations, 
and from this perspective they are exactly their geometric counterparts. (For a quick 
introduction to the Lie theoretic background, see e.g. \cite{prehom}).
In type $A_n$, we get the usual Grassmannians:

\setlength{\unitlength}{4mm}
\thicklines
\begin{picture}(18,3)(-1,.9)
\multiput(0,2)(2,0){7}{$\circ$}
\put(4,2){$\bullet$}
\multiput(0.3,2.2)(2,0){6}{\line(1,0){1.8}}
\put(14.2,2){$\simeq$}
\put(16,2){$G/P=G(3,8)\subset\PP(\wedge^3\CC^8)$}
\end{picture}

In types $B_n, D_n$ (resp. $C_n$) there is an invariant quadratic (resp. symplectic) form
preserved by the group $G$, 
and the generalized Grassmannians $G/P_k$ are $OG(k,m)$, for $m=2n+1$ or $2n$ 
(resp. $IG(k,2n)$), the subvarieties of the usual Grassmannians parametrizing 
{\it isotropic subspaces}. 

This has to be taken with a grain of salt for $k=n$ or 
$n-1$ in type $D_n$ (and also for $k=n$ in type $B_n$): the variety of isotropic spaces
$OG(n,2n)\subset\PP(\wedge^n\CC^{2n})$ has two connected components; moreover, the restriction
$L$ of the Pl\"ucker line bundle to a component is divisible by two, and the line bundle
$M$ such that $L=2M$ embeds this component into the projectivization of a {\it half-
spin representation} $\Delta_{2n}$. As embedded varieties, these two components are in fact  undistinguishable.
We will denote them by $\SS_{2n}\subset \PP(\Delta_{2n})$ and call them the {\it spinor 
varieties}. 
 
\setlength{\unitlength}{4mm}
\thicklines
\begin{picture}(18,4.2)(-5,.9)
\multiput(2,2)(2,0){4}{$\circ$}
\multiput(2.35,2.2)(2,0){3}{\line(1,0){1.75}}
\put(9.8,.25){$\circ$}
\put(9.8,3.75){$\bullet$}
\put(8.4,2.35){\line(1,1){1.5}}
\put(8.4,2.05){\line(1,-1){1.5}}
\put(11.4,3.6){$\simeq$}
\put(13.1,3.6){$G/P=\SS_{12}\subset\PP(\Delta_{12})$}
\end{picture}

\subsection{Some reasons to care}
Being homogeneous, homogeneous spaces could look boring! But we are far from knowing everything
about them. For example: 
\begin{enumerate}
\item Their Chow rings are not completely understood, except for usual Grassmannians and a 
few other varieties: {\it Schubert calculus} has been an intense field of research 
since the 19th century,
involving geometers, representation theoretists, and combinatorists! Modern versions include
K-theory, equivariant cohomology, quantum cohomology, equivariant quantum K-theory... Interesting 
{\it quantum to classical} principles are known in some cases, which for instance allow to 
deduce the quantum cohomology of Grassmannians from the usual intersection theory on two-step
flag varieties. 
\item Derived categories are fully described only for special cases (Grassmannians,
quadrics, isotropic Grassmannians of lines, a few other sporadic cases), although 
important progress have been made for classical Grassmannians. 
\item Characterizations of homogeneous spaces are important but not known (except in small dimension). This is one of the potential interests of rational homogeneous spaces for the 
algebraic geometer: their behavior is often extremal in some sense. Here are some important conjectures:
\end{enumerate}

\smallskip\noindent {\bf Campana-Peternell conjecture (1991)}. {\it Let $X$ be a smooth 
complex Fano variety with nef tangent bundle. Then $X$ is homogeneous}. 

\medskip
See \cite{cp} for a survey. Nefness is a weak version of global generation: varieties
with globally generated tangent bundles are certainly homogeneous. Note that varieties
with ample tangent bundles are projective spaces: this was conjectured by Hartshorne
and Frankel, and proved by Mori.

\smallskip
A smooth codimension one distribution on a variety $Y$ is defined as a corank one sub-bundle 
$H$ of the tangent bundle $TY$. Let $L=TY/H$ denote the quotient line bundle. The Lie bracket
on $TY$ induces a linear map $\wedge^2H\lra L$. This gives what is called a {\it contact structure}
when this skew-symmetric form is non degenerate (this implies that the dimension $d=2n+1$ of
$Y$ is odd and that the canonical line bundle $\omega_Y=-(n+1)L$, see \cite{kpsw}) . 

\smallskip\noindent {\bf Lebrun-Salamon conjecture (1994)}. {\it Let $Y$ be a smooth complex Fano
variety admitting a contact structure. Then $Y$ is homogeneous}.

\medskip Hartshorne's conjecture (1974) 
in complex projective geometry is that if $Z\subset \PP^N$ 
is smooth, non linearly degenerate, of dimension $n>\frac{2N}{3}$, then $Z$ must be a 
complete intersection. In particular it must be linearly normal, meaning that it is not a 
projection from $\PP^{N+1}$. Although Hartshorne's conjecture is still widely open, 
this last statement was indeed shown by Zak under the weaker assumption that $n>\frac{2N-2}{3}$. 
He also proved \cite{zak}: 

\smallskip\noindent {\bf Zak's theorem II (1981)}. {\it Let $Z$ be a smooth, non linearly degenerate, non linearly normal subvariety of $\PP^N$, of dimension $n=\frac{2N-2}{3}$. Then 
$Z$ is homogeneous. }

Those special homogeneous spaces are called Severi varieties (see section \S 1.5). 

\subsection{Homogeneous spaces and Fano manifolds}

The importance of homogeneous spaces in the theory of Fano varieties comes from the following 

\smallskip\noindent {\bf Fact}. \emph{Generalized Grassmannians are Fano manifolds of high index.}  

\begin{definition}
Let $X$ be a Fano manifold of dimension $n$ and Picard number one, so that $Pic(X)=\ZZ L$ 
for some ample line bundle $L$, and $K_X=-i_X L$. Then $i_X$ is the index and the coindex 
$c_X:= n+1-i_X$.
\end{definition} 

\begin{theorem} $\mathrm{(Kobayashi-Ochiai, 1973).}$

\begin{enumerate}
\item If $c_X=0$, then $X\simeq\PP^n$ is a projective space.
\item If $c_X=1$, then $X\simeq\QQ^n$ is a quadric.
\end{enumerate}
\end{theorem}

Here is a list of generalized Grassmannians of coindex $2$ (del Pezzo manifolds) and $3$ (Mukai varieties). 
\begin{center}
$$\begin{array}{cccc} 
 & & \mathrm{dimension} & \mathrm{index} \\
del\; Pezzo & G(2,5) & 6 & 5 \\
Mukai & G_2/P_2& 5 & 3 \\
       & IG(3,6) & 6 & 4 \\
        & G(2,6) & 8 & 6 \\
        & \SS_{10} & 10 & 8 \\ 
\end{array}$$
\end{center}

\medskip \noindent 
{\bf Observation}. Suppose $L$ is very ample and embeds $X$ in $\PP V$. Then a 
smooth hyperplane section $Y=X\cap H$ has the same coindex $c_Y=c_X$. 

\smallskip As a consequence, a generalized Grassmannian contains, as linear sections, a large
family of Fano submanifolds of dimension down to $c_X$. Moreover, these families are always 
{\it locally complete}, in the sense that every small complex deformation of any member of 
the family is of the same type. In particular, Mukai varieties allow to construct locally 
complete families of prime Fano threefolds of index $1$.

\begin{theorem}[Fano-Iskhovskih's classification \cite{ip}]
Any prime Fano threefold of index one is either:
\begin{enumerate}
\item a complete intersection in a weighted projective space,
\item a quadric hypersurface in a del Pezzo homogeneous space,
\item a linear section of a homogeneous Mukai variety,
\item a trisymplectic Grassmannian.
\end{enumerate}
\end{theorem}

For the last case, one gets Fano threefolds in $G(3,7)$, parametrizing  spaces that are isotropic with respect to three alternating two-forms. 
A unified description of the last three cases is by zero loci of 
sections of homogeneous bundles on homogeneous spaces; this stresses
the importance of understanding those bundles.

\subsection{Correspondences and lines}

There are natural projections from (generalized) flag manifolds to Grassmannians. 
In particular, given two generalized Grassmannians $G/P$ and $G/Q$ of the same complex
Lie group $G$, we can connect them through a flag manifold $G/(P\cap Q)$ (we may suppse that 
$P\cap Q$ is still parabolic):
$$\xymatrix{
 &  G/(P\cap Q) \ar@{->}[ld]_{p}\ar@{->}[rd]^{q}  &  \\
 G/P & & G/Q 
}$$
The fiber of $p$ is $P/(P\cap Q)$, in particular it is homogeneous. If its automorphism group is
the semisimple Lie group $L$, we can write $P/(P\cap Q)=L/R$ and this generalized flag manifold
is determined as follows. 
\begin{itemize}
\item Start with the Dynkin diagram $\Delta$ of $G$, with the two nodes $\delta_P$ and $\delta_Q$ 
defining $G/P$ and $G/Q$. 
\item Suppress the node $\delta_P$ and the edges attached to it. 
\item Keep the connected component of the remaining diagram that contains $\delta_Q$. 
\end{itemize}
The resulting marked Dynkin diagram is that of the fiber $P/(P\cap Q)=L/R$. 

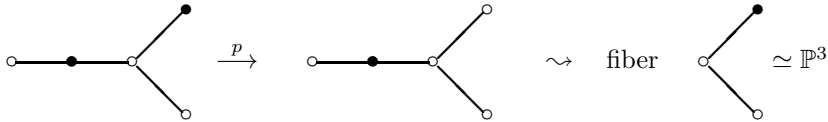
\setlength{\unitlength}{4mm}
\thicklines
\begin{picture}(18,4.2)(3,.9)
\multiput(4,2)(2,0){3}{$\circ$}
\multiput(4.35,2.2)(2,0){2}{\line(1,0){1.75}}
\put(9.8,.25){$\circ$}
\put(9.8,3.75){$\bullet$}\put(6,2){$\bullet$}
\put(8.4,2.35){\line(1,1){1.5}}
\put(8.4,2.05){\line(1,-1){1.5}}
\put(11,2){$\stackrel{p}{\longrightarrow}$}
\multiput(14,2)(2,0){3}{$\circ$}
\multiput(14.35,2.2)(2,0){2}{\line(1,0){1.75}}
\put(19.8,.25){$\circ$}
\put(16,2){$\bullet$}\put(19.8,3.75){$\circ$}
\put(18.4,2.35){\line(1,1){1.5}}
\put(18.4,2.05){\line(1,-1){1.5}}
\put(22,2){$\leadsto$}
\put(24,2){fiber}
\put(27.4,2.35){\line(1,1){1.5}}
\put(27,2){$\circ$}
\put(27.4,2.05){\line(1,-1){1.5}}
\put(28.8,3.75){$\bullet$}
\put(28.8,.25){$\circ$}
\put(29.5,2){$\simeq\PP^3$}
\end{picture}
\smallskip

This rule has an obvious extension to the case where $G/P$ and $G/Q$ are not necessarily
generalized Grassmannians. 

In particular, start with a generalized Grassmannian $G/P$, defined 
by the marked Dynkin diagram $(\Delta,\delta_P)$. Let $\delta_P^{prox}$ be the set of vertices in $\Delta$ that are connected to $\delta_P$. Let $G/P_{prox}$ be the generalized flag manifold
defined by the marked Dynkin diagram $(\Delta,\delta_P^{prox})$. Then the fibers of $q$ are
projective lines! 

\begin{theorem}\cite{LMproj}
If $\delta_P$ is a long node of $\Delta$ (in particular, if $\Delta$ is simply laced),
then $G/P_{prox}$ is the variety of lines on $G/P$. 

If $\delta_P$ is a short node of $\Delta$, then the variety of lines on $G/P$ is irreducible
with two $G$-orbits, and $G/P_{prox}$ is the closed one. 
\end{theorem}

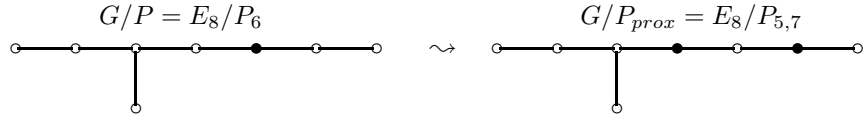
\setlength{\unitlength}{4mm}
\thicklines
\begin{picture}(18,5.9)(-1,-1.5)
\multiput(0,2)(2,0){7}{$\circ$}\put(4,0){$\circ$}
\put(8,2){$\bullet$}
\multiput(0.3,2.2)(2,0){6}{\line(1,0){1.8}}
\put(4.2,0.3){\line(0,1){1.8}}
\put(3,3){$G/P=E_8/P_6$}
\put(14,2){$\leadsto$}
\multiput(16,2)(2,0){7}{$\circ$}\put(20,0){$\circ$}
\put(22,2){$\bullet$}\put(26,2){$\bullet$}
\multiput(16.3,2.2)(2,0){6}{\line(1,0){1.8}}
\put(20.2,0.3){\line(0,1){1.8}}
\put(19,3){$G/P_{prox}=E_8/P_{5,7}$}
\end{picture}

In particular, any generalized Grassmannians parametrizes linear spaces on a generalized Grassmannian defined by an extremal node of a Dynkin diagram. 
\smallskip

There is a similar statement for the variety of lines passing through a fixed point $x$
of $G/P\subset \PP V$. These lines are parametrized by a subvariety $\Sigma_x\subset\PP(T_xX)$,
independant of $x$ up to projective equivalence. Moreover $\Sigma_x$ is stable under the action
of the stabilizer $P_x\simeq P$ of $x$. The semisimple part $H$ of $P$ has Dynkin diagram $\Delta_H=
\Delta - \{\delta_P\}$, and we can define a parabolic subgroup $Q\subset H$ by marking in 
$\Delta_H$ the vertices that  in $\Delta$ were connected to $\delta_P$.

\begin{theorem}\cite{LMproj}
If $\delta_P$ is a long node of $\Delta$ (in particular, if $\Delta$ is simply laced),
then $\Sigma_x\simeq H/Q$.
If $\delta_P$ is a short node of $\Delta$, then $\Sigma_x$ is made of two  $H$-orbits, and $H/Q$ is the closed one. 
\end{theorem}

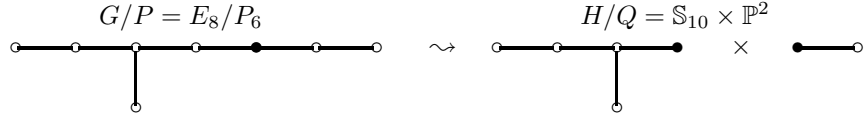
\setlength{\unitlength}{4mm}
\thicklines
\begin{picture}(18,5.9)(-1,-1.5)
\multiput(0,2)(2,0){7}{$\circ$}\put(4,0){$\circ$}
\put(8,2){$\bullet$}
\multiput(0.3,2.2)(2,0){6}{\line(1,0){1.8}}
\put(4.2,0.3){\line(0,1){1.8}}
\put(3,3){$G/P=E_8/P_6$}
\put(14,2){$\leadsto$}
\multiput(16,2)(2,0){3}{$\circ$}\put(20,0){$\circ$}\put(28,2){$\circ$}
\put(22,2){$\bullet$}\put(26,2){$\bullet$}\put(24,2){$\times$}
\multiput(16.3,2.2)(2,0){3}{\line(1,0){1.8}}\put(26.3,2.2){\line(1,0){1.8}}
\put(20.2,0.3){\line(0,1){1.8}}
\put(19,3){$H/Q=\SS_{10}\times\PP^2$}
\end{picture}

\smallskip
Variety of lines through fixed points in homogeneous spaces are 
instances of the so-called VMRTs (varieties of minimal rational tangents), which
have been extensively studied in the recent years in connexion with 
rigidity problems (see e.g. \cite{hwang} for an introduction). 

From deformation theory, one knows that the variety of lines passing through a general point
of a Fano manifold $X$ has dimension $i_X-2$. With the previous notations, this gives a nice
way to compute the index of a generalized Grassmannian (at least for the long node case):
 $$ i_{G/P}=\dim (H/Q)+2.$$

\noindent{\it Example 1}. Let us consider the rank two (connected) Dynkin diagrams. 
\begin{enumerate} 
\item $\Delta=A_2$, $G=PGL_3$. Then the two generalized Grassmannians are projective planes,
and each of them parametrizes the lines in the other one. 
\item $\Delta=B_2$, $G=PSO_5$. The two generalized Grassmannians are the quadric $\QQ^3$
and $OG(2,5)$. The latter parametrizes the lines in $\QQ^3$. But note that $OG(2,5)$ has dimension $3$ and anticanonical twice the restriction of the Pl\"ucker line bundle, hence index $4$ since 
this restriction is $2$-divisible. So by the Kobayashi-Ochiai theorem $OG(2,5)\simeq\PP^3$! Moreover, lines
in $\PP^3$ are parametrized by $G(2,4)\simeq\QQ^4$, and $\QQ^3$ is the closed $PSO_5$-orbit. 
\item $\Delta=C_2$, $G=PSp_4$. The two generalized Grassmannians are $\PP^3$
and $IG(2,4)$. The latter is a hyperplane section of $G(2,4)\simeq\QQ^4$, hence a copy of $\QQ^3$.
Of course we recover the previous case.
\item $\Delta=G_2$, $G=G_2$. The two generalized Grassmannians $G_2/P_1$ and $G_2/P_2$ have
the same dimension $5$, but different indexes $5$ and $3$. In particular, by the Kobayashi-Ochiai
theorem
$G_2/P_1\simeq\QQ^5$. Moreover $G_2/P_2$ must be the closed $G_2$-orbit in the variety of lines
in $\QQ^5$, which is $OG(2,7)$. 
\end{enumerate}

\subsection{Sporadic examples}
There exist a few series of generalized Grassmannians with strikingly similar properties. 

\medskip\noindent   {\it Severi varieties}. These are the four varieties
$$v_2(\PP^2)\subset \PP^5, \quad \PP^2 \times\PP^2\subset \PP^8, \quad G(2,6)\subset  \PP^{14}, 
 \quad E_6/P_1\subset  \PP^{26}.$$
Here $v_2$ means that we consider the second Veronese embedding.
Each of these varieties is the singular locus of a special cubic hypersurface (the secant 
variety), and the derivatives of this cubic define a quadro-quadric Cremona tranformation. 
Note also that they are varieties of dimension $2a$ embedded inside $\PP^{3a+2}$ for 
$a=1,2,4,8$. These are the homogeneous spaces that appear in Zak's Theorem II \cite{zak}.

\medskip\noindent   {\it Legendrian varieties}. These are the four varieties
$$IG(3,6)\subset \PP^{13}, \quad G(3,6)\subset \PP^{19}, \quad \SS_{12}\subset  \PP^{31}, 
 \quad E_7/P_1\subset  \PP^{55}.$$
Each of these varieties is contained the singular locus of a special quartic hypersurface 
(the tangent variety), and the derivatives of this quartic define a cubo-cubic Cremona tranformation. Moreover they have the remarkable {\it one apparent double point property}, 
which means that through a general point of the ambient projective space passes exactly 
one bisecant. 
Note also that they are varieties of dimension $3a+3$ embedded inside $\PP^{6a+9}$ for 
$a=1,2,4,8$. Finally, their varieties of lines through a given point are 
nothing else than the Severi varieties
\cite{LMleg}!

\section{Borel-Weil theory and applications} 
Let $G/P$ be a generalized Grassmannian and suppose that $G$
is simply connected. Let $L$ be the ample generator of the Picard group. Then $L$ is in fact
very ample and $G$-linearizable. In particular $H^0(G/P,L)$ is a $G$-module and $G/P$ is
$G$-equivariantly embedded inside the dual linear system $|L|^\vee = \PP H^0(G/P,L)^\vee$.
A typical example is the Pl\"ucker embedding of a Grassmannian (see e.g. \cite{fh} for basics
in representation theory).

\subsection{The Borel-Weil theorem}
More generally, on any generalized flag manifold $G/Q$, any line bundle $M$ is
$G$-linearizable, so $H^0(G/Q,M)$ is a $G$-module. Moreover $M$ is generated by global sections as soon as it is nef. 

\begin{theorem}[Borel-Weil]
For any nef line bundle $M$ on $G/Q$, 
the space of sections $H^0(G/Q,M)$ is an irreducible $G$-module.
\end{theorem}

A line bundle $M$ on $G/Q$ is defined by a weighted version of the marked Dynkin 
diagram that defines $G/Q$. Moreover, it is nef/globally generated (resp. ample/very ample) 
exactly when the weights are non negative (resp. positive), and then the $G$-module 
$H^0(G/Q,M)$ is defined by the {\it same} weighted diagram. 

\smallskip
Starting from a projection $p: G/(P\cap Q)\lra G/P$, by homogeneity the sheaf $E_P=p_*M$ is a 
$G$-equivariant vector bundle on $G/P$.  Symmetrically 
there is a vector bundle $E_Q=q_*M$ on $G/Q$, and  
$$H^0(G/Q,E_Q)=H^0(G/(P\cap Q),M)=H^0(G/P,E_P).$$
Such identifications allow to play with sections in different homogeneous 
spaces and describe nice 
correspondences between their zero loci. 

\medskip\noindent{\it Example 2}. Consider the flag manifold $F(2,3,5)$ with its two projections $p$ and $q$ to 
$G(2,5)$ and $G(3,5)$. The tautological and quotient bundles $U_2, Q_2$ on $G(2,5)$, 
$U_3, Q_3$ on $G(3,5)$, pull-back to vector bundles on $F(2,3,5)$ for which we keep the 
same notations. The minimal very ample line bundle on $F(2,3,5)$ is $L=\det(U_2)^\vee
\otimes \det(Q_3)$. Its push-forwards to the two Grassmannians are 
$E_2=\det(U_2)^\vee\otimes\wedge^2Q_2\simeq Q_2^*(2)$ on $G(2,5)$ and 
$E_3=\wedge^2U_3^\vee\otimes\det(Q_3)\simeq U_3(2)$ on $G(3,5)$. 
These are two vector bundles with determinant $\cO(5)$. As a consequence, a general
section $s$ of $L$ defines two Calabi-Yau threefolds
$$Z_2(s)\subset G(2,5) \qquad \mathrm{and} \qquad  Z_3(s)\subset G(3,5).$$

\begin{proposition}\cite{kr} 
The Calabi-Yau threefolds $Z_2(s)$ and $Z_3(s)$ are derived equivalent, but not  birationally equivalent in general.
\end{proposition}

Derived equivalent means that their derived categories of coherent sheaves are  equivalent
as triangulated categories. For two smooth projective varieties $X_1$ and $X_2$, this is 
a very strong property, which implies that they are in fact isomorphic as soon as one of them has 
an ample or anti-ample canonical bundle (Bondal-Orlov). Non isomorphic but derived equivalent 
varieties with trivial canonical bundle recently attracted considerable attention, and the 
previous example of Fourier-Mukai partners is among the simplest. 

\medskip\noindent{\it Example 3}. 
A three-form $\omega\in \wedge^3(\CC^n)^\vee$ defines a 
section of $\wedge^3T^\vee$ on $G(k,n)$ for any $k\ge 3$, and sections 
of $Q^\vee(1)$ on $G(2,n)$ and $\wedge^2Q^\vee(1)$ on $\PP^{n-1}$. 
The latter leads to Pfaffian loci in $\PP^{n-1}$. The previous one gives 
{\it congruences of lines} in $G(2,n)$: the zero locus of a general section
of $Q^\vee(1)$ is a $(n-2)$-dimensional prime (for $n\neq 6$) Fano manifold of index $3$. For $n=6$, one actually gets $\PP^2\times \PP^2$. 

\smallskip
For $n=7$, this congruence of lines is isomorphic with $G_2/P_2$. Indeed, 
the stabilizer in $SL_7$ of a general $\omega\in \wedge^3(\CC^7)^\vee$ is isomorphic to $G_2$.
This stabilizer acts on the associated congruence, and since there is no non trivial 
closed $G_2$-orbit of dimension smaller than five, this congruence has to be a $G_2$-Grassmannian. But it is not $G_2/P_1=\QQ^5$, whose index
is $5$, so it must be $G_2/P_2$.

\smallskip
For $n=8$, a general $\omega\in \wedge^3(\CC^8)^\vee$ is equivalent to the three-form $$\omega_0(x,y,z)=\mathrm{trace}(x[y,z])$$ on $\fsl_3$, which is stabilized by the adjoint action of $SL_3$. 
The congruence of lines defined by $\omega_0$ is the variety of {\it abelian} planes inside $\fsl_3$, a smooth
compactification of $SL_3/T$, for $T\subset SL_3$ a maximal torus: so a quasi-homogeneous, but not homogeneous sixfold.

\medskip These constructions intend to illustrate the general 

\smallskip\noindent {\bf Principle}. \emph{Vector bundles on generalized flag manifolds allow to easily construct interesting varieties, notably Fano or Calabi-Yau.}


\subsection{Flops}
Consider a rank two Dynkin diagram $\Delta$, with corresponding group $G$. The flag varieties 
of $G$ are the complete flag variety $G/B$ and the two generalized Grassmannians $G/P_1$ and $G/P_2$. Let $L$ denote the minimal ample line bundle on $G/B$. Its push-forwards to $G/P_1$ and $G/P_2$
are two vector bundles $E_1$ and $E_2$.

$$\xymatrix{& & L^\vee\ar@{->}[d] & &\\
E_1^\vee\ar@{->}[d] & & G/B \ar@{->}[lld]_{p_1}\ar@{->}[rrd]^{p_2}  & & E_2^\vee\ar@{->}[d] \\
 G/P_1 &&  && G/P_2 
}$$

Consider a point in the total space of $L^\vee$, that is a pair $(x,\phi)$ with 
$x\in G/B$ and $\phi$ a linear form on the fiber $L_x$. Let $y=p_1(x)$ and $e\in E_{1,y}$, the 
corresponding fiber of $E_1$. Since $E_1=p_{1*}L$, we can see $e$ as a section of $L$ over  $p_1^{-1}(y)$, then evaluate this section at $x$, and apply $\phi$. This defines a morphism 
$v_1: Tot(L^\vee)\lra Tot(E_1^\vee)$. Moreover two points $(x,\phi)$ and $(x',\phi')$ have
the same image if and only if $p_1(x)=p_1(x')$ and $\phi=\phi'=0$. This yields:

\begin{proposition} The morphism $v_i: Tot(L^\vee)\lra Tot(E_i^\vee)$ is the blowup of $G/P_i$. \end{proposition}

In particular we get a nice birational map between $Tot(E_1^\vee)$ and $Tot(E_2^\vee)$. 
In type $A_1\times A_1$ this is the classical {\it Atiyah flop} (in dimension three), 
and in type $A_2$ the classical {\it Mukai flop} (dimension four). In type $B_2=C_2$ 
this is the {\it Abuaf-Segal flop} (dimension five) \cite{segal}, and in type $G_2$ the {\it Abuaf-Ueda flop} (dimension seven) \cite{ueda}.  An unusual feature of the two latter flops is 
that they are quite non symmetric, the exceptional loci on both sides being rather different. 

\begin{theorem}
All these flops are derived equivalences.
\end{theorem}

This supports a famous conjecture of Bondal and Orlov according to which varieties connected 
by flops should always be derived equivalent.

\medskip 
The three bundles $L, E_1, E_2$ have the same space of global sections $V$. So if we pick a (general) section $s\in V$, we get three zero-loci $Z(s), Z_1(s), Z_2(s)$ and a diagram:
$$\xymatrix{& & L\ar@{->}[d] & &\\
E_1\ar@{->}[d] & & G/B \ar@{->}[lld]_{p_1}\ar@{->}[rrd]^{p_2}  & & E_2\ar@{->}[d] \\
 G/P_1 && Z(s)\ar@{^{(}->}[u]\ar@{->}[ll]_{u_1}\ar@{->}[rr]^{u_2} && G/P_2 \\ 
 Z_1(s)\ar@{^{(}->}[u] &&&& Z_2(s)\ar@{^{(}->}[u]
}$$
Observe that $p_1$ and $p_2$ are $\PP^1$-fibrations, such that $L$ is a relative hyperplane
bundle. So $E_1$ and $E_2$ are rank two bundles. Moreover, $s$ vanishes either at a unique point 
of a fiber, or on this whole fiber. This implies that $u_1$ and $u_2$ are birational morphisms,
with exceptional divisors that are $\PP^1$-bundles over the codimension two subvarieties $Z_1(s)$
and $Z_2(s)$. More precisely:

\begin{proposition}\label{blowup} 
The morphism $u_i: Z(s)\lra G/P_i$ is the blowup of $Z_i(s)$. \end{proposition}

Note that the canonical bundle of $G/B$ is $-2L$, so that $Z(s)$ is Fano: this is 
one of the (not so many) examples of blowup of a Fano manifold that remains Fano. 

\smallskip
In type $A_2$, $Z(s)$ is a del Pezzo surface of degree $6$, each $Z_i(s)$ consists in three
points in a projective plane, and the birationality between the two planes is the classical
Cremona transformation. In type $C_2$, $Z_1(s)$ is a quintic elliptic curve $E$ 
in $C_2/P_1=\PP^3$, and we obtain one of the examples of the Mori-Mukai classification 
of Fano threefolds with $b_2=2$. The birational map to $C_2/P_2=\QQ^3$ is defined by the linear system of cubics through $E$. Moreover $Z_2(s)$ is again an elliptic curve, isomorphic 
to $E$ since both curves can be identified with the intermediate Jacobian of the Fano threefold $Z(s)$. More to come about type $G_2$.

\medskip\noindent {\it Remark}.
In general, consider a projection 
$p: G/P\lra G/Q$ and a line bundle $L$ on $G/P$. The vector bundle $E=p_*L$ is a
$G$-equivariant bundle, defined by a $P$-module $W$. 
Consider a section $s$ of $L$, vanishing along a smooth hypersurface $Z(s)\subset G/P$. The fiber over $x$ of the projection $p_s : Z(s)\lra G/Q$ is the zero locus of the restriction $s_x$ of $s$ to the fibers of $G/P$. Its isomorphism type only depends on the $P$-orbit of $s_x$ in $W$. This is
one of the many reasons to be interested in representations with finitely many orbits.

\subsection{Consequences for the Grothendieck ring of varieties}

The Grothendieck ring of (complex) varieties is the ring generated by isomorphism classes
of algebraic varieties (no scheme structures) modulo the relations
\begin{enumerate}
\item $[X]=[Y]+[X\backslash Y]$ for $Y\subset X$ a closed subvariety,
\item $[X\times X']=[X]\times [X']$.
\end{enumerate}
An easy consequence is that if $X$ is a Zariski locally trivial fiber bundle over $Z$,
with fiber $F$, then $[X]=[Z]\times [F]$. 
If we denote by $\LL$ the class of the affine line, the usual cell 
decomposition of projective spaces yields
$$[\PP^n]=1+\LL+\cdots +\LL^n=\frac{1-\LL^{n+1}}{1-\LL}.$$
Since the complete flag manifold $\FF_n=SL_n/B$ is a composition of projective bundles, one 
deduces that 
$$[\FF_n]=\frac{(1-\LL^2)\cdots (1-\LL^{n})}{(1-\LL)^{n-1}}.$$
Similar formulas exist for any $G/P$: they admit stratifications by affine spaces (cell decompositions) and the
classes $[G/P]$ are therefore polynomials in $\LL$, which admit nice factorizations as rational
functions. 

\smallskip
Let us come back to Proposition \ref{blowup}. Since the the blowup $u_i$ gives a $\PP^1$-bundle over $Z_i(s)$, and an isomorphim over its complement,  we get that 
$$[Z(s)]=[G/P_i-Z_i(s)]+[Z_i(s)]\times [\PP^1]=[G/P_i]+[Z_i(s)]\times \LL.$$
We deduce the identity 
$$([Z_1(s)]-[Z_2(s)])\times \LL = [G/P_2]-[G/P_1]=0.$$
(In fact $G/P_1$ and $G/P_2$ have the same Hodge numbers as a projective space: such 
varieties are called minifolds). This has the unexpected consequence
that \cite{imou}:

\begin{theorem}\label{zerodiv} 
$\LL$ is a zero divisor in the Grothendieck ring of varieties.
\end{theorem}

\proof Let $G=G_2$. One needs to check that in this case, $[Z_1(s)]-[Z_2(s)]\ne 0$. This 
requires more sophisticated arguments.  First, one shows that for $s$ general, $Z_1(s)$ 
and $Z_2(s)$ are smooth Calabi-Yau threefolds \cite{imou}. 
Second, their Picard groups are both 
cyclic, but the minimal ample generators have different degrees. So $Z_1(s)$ and $Z_2(s)$
are not isomorphic. 

Therefore they are not birational, because Calabi-Yau's are minimal
models, so they would be isomorphic in codimension one, and since the Picard groups are
cyclic their minimal generators would match. Therefore they are not stably birational 
either, because if $Z_1(s)\times\PP^m$ was birational to $Z_2(s)\times\PP^m$, then their MRC fibrations, or maximal rationally connected fibrations, would also be birational; but since 
they are not uniruled, $Z_1(s)$ and $Z_2(s)$ are the basis of these fibrations. 
Finally, a deep result of Larsen and Lunts implies that $[Z_1(s)]-[Z_2(s)]$
does not belong to the ideal generated by $\LL$, and in particular it must be non zero. \qed

\medskip The fact that the Grothendieck ring is not a domain was first observed by Poonen 
(2002). That $\LL$ is a zero divisor was first shown by Borisov using the Pfaffian-Grassmannian
equivalence \cite{bo} (see Proposition \ref{bc}). 

\subsection{Representations with finitely many orbits} Consider a semisimple Lie group $G$, and 
an irreducible representation $W$ such that $\PP (W)$ contains only finitely many $G$-orbits. Over $\CC$ these representations 
were classified by Kac \cite{Kac80}, who proved that most of them are {\it parabolic}. Parabolic representations
are exactly those representations spanned by the varieties of lines through a given point of 
a generalized Grassmannian. These representations always admit finitely many orbits. 
This implies that there are  very strong connections between rational
homogeneous spaces and prehomogeneous spaces \cite{prehom}.

\setlength{\unitlength}{4mm}
\thicklines
\begin{picture}(18,5)(-1,-.5)
\multiput(0,2)(2,0){7}{$\circ$}\put(4,0){$\circ$}
\put(8,2){$\bullet$}
\multiput(0.3,2.2)(2,0){6}{\line(1,0){1.8}}
\put(4.2,0.3){\line(0,1){1.8}}
\put(3,3){$G/P=E_8/P_6$}
\put(15,2){$\leadsto$}
\put(18,2){$\Delta_{16}\times\CC^3\curvearrowleft Spin_{10}\times GL_3$}
\end{picture}

\smallskip\noindent {\it Example 4}. For dimensional reasons the spaces of three-forms 
$\wedge^3(\CC^n)^\vee$ cannot have finitely many $GL_n$-orbits
when $n\ge 9$. They do have finitely many $GL_n$-orbits for $n\le 8$ because they are parabolic,
coming from generalized Grassmannians of type $E_n$. Orbits were classified long ago.

Note that each orbit is a locally closed subvariety, whose boundary is a union of smaller
dimensional orbits. So there is a natural Hasse diagram encoding, for each orbit closure, 
which other ones  are the components of its boundary. 

\medskip
For $n=6$ the orbit closures are particularly simple to describe. 
The Hasse diagram is just a line
$$\xymatrix{  Y_0\ar@{->}[r]&  Y_1\ar@{->}[r]&  Y_5\ar@{->}[r]& Y_{10}\ar@{->}[r] & Y_{20},}$$ 
where the index is the codimension:

$$\begin{array}{lcl}
 Y_0 & \quad & ambient\; space \\
 Y_1 & & tangent\; quartic\; hypersurface \\
 Y_5 & & singular\: locus\; of\;the\; quartic \\
 Y_{10} && cone\; over\; the\; Grassmannian\; G(3,6)\\
 Y_{20} && origin
 \end{array}$$

Alternatively, $Y_5$ can be described as parametrizing those three-forms $\theta$ that factor 
as $\omega\wedge\ell$ for some vector $\ell$ and some two-form $\omega$. In general the 
line $L=\langle\ell\rangle$ is uniquely defined by $\theta$. Moreover $\omega$ (which is defined modulo $L$) is a two-form of maximal rank $4$, supported on a uniquely defined hyperplane $H$. 
Then $\theta$ belongs to $\wedge^2H\wedge L$. We conclude that $Y_5$ is birational to the total
space of the rank six vector bundle $E_0=\wedge^2U_5\wedge U_1$ over the flag variety $F(1,5,6)$. 

But we can also forget $U_1$ and $U_5$ and  we get a diagram 
 $$
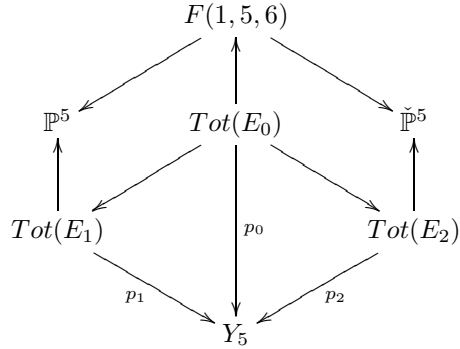
\xymatrix{&  F(1,5,6) \ar@{->}[dl]\ar@{->}[dr]& \\
 \PP^5 & Tot(E_0)\ar@{->}[ld]\ar@{->}[u]\ar@{->}[dd]^{p_0}\ar@{->}[rd] & \check{\PP}^5 \\
 Tot(E_1)\ar@{->}[rd]_{p_1}\ar@{->}[u] && Tot(E_2)\ar@{->}[ld]^{p_2}\ar@{->}[u] \\
 & Y_5 &
}$$
where $E_1$ is the vector bundle $\wedge^2\CC^6\wedge U_1$ and $E_2$ is the vector bundle 
$\wedge^3U_5$, both of rank ten. It is a nice exercise to describe this flop. One checks in 
particular that

\begin{proposition}
$p_0,p_1,p_2$ are resolutions of singularities. $p_0$ is divisorial, while $p_1$ and $p_2$ are crepant resolutions, with exceptional loci of codimension three. 
\end{proposition}

\medskip For $n=7$ the situation is  slightly more complicated. Here is the Hasse diagram, where as
previously the index stands for the codimension. 
$$
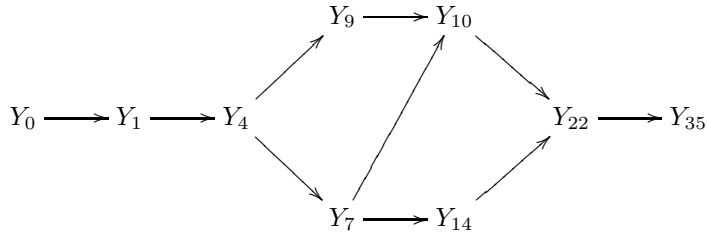
\xymatrix{ &&& Y_9\ar@{->}[r] & Y_{10} \ar@{->}[dr] & & \\
Y_0 \ar@{->}[r] &  Y_1 \ar@{->}[r]  & Y_4  \ar@{->}[rd]\ar@{->}[ru] 
&& &  Y_{22} \ar@{->}[r] & Y_{35} \\
 &&& 
 Y_7 \ar@{->}[r] \ar@{->}[uur] & Y_{14} \ar@{->}[ur] &&&
 }$$
 
For $n=8$ there are $22$ orbits, partially described in \cite{coble}.

\subsection{Applications to linear sections} 
A nice case is when the generalized Grassmannian $G/P$ is defined by a node that 
splits the Dynkin diagram $\Delta$ into the union of two diagrams, one of which is 
of type $A_{k-1}$, with $k\ge 1$ ($k=1$ occurs when the node is at an end of the diagram). 
This means that the lines through some fixed point in $G/P$ are 
parametrized by $\PP^{k-1}\times S$ for some homogeneous space $S$ (of some smaller 
Lie group $H$). 

If the minimal homogeneous embedding of $S$ is inside $\PP (U)$ for 
some irreducible $H$-module $U$, the minimal homogeneous embedding of $\PP^{k-1}\times S$
is inside $\PP (\CC^{k}\otimes U)$. Moreover, by Kac's results mentionned above, the action of 
$GL_{k}\times H$ on $\CC^{k}\otimes U$ admits finitely many orbits. This implies 
that the Grassmannian $G(k,U)$ also admits finitely many $H$-orbits, and therefore:

\begin{proposition} 
There exist only finitely many isomorphism types of codimension $k$ linear sections 
of $S\subset\PP (U)$. 
\end{proposition}

Of course the open orbit yields smooth linear sections, but it may happen that other orbits
also give smooth sections, although degenerate in some sense. A very interesting case is that 
of the spinor variety $\SS_{10}$ \cite{kuzspin,bfm}. 
Starting from the Dynkin diagram of type $E_8$, we see that
linear sections of codimension up to three have finitely many isomorphism types. One can check 
that there exists only one type of smooth sections of codimension $1$, 
but two different types of codimension 
$2$ and four types of codimension $3$. 

\smallskip
This gives examples of non locally rigid Fano varieties of high index. Conversely, the 
local deformations of Fano linear sections of generalized Grassmannians are always 
of the same type, and the generic Fano linear sections are locally rigid if and only 
if their marked Dynkin diagrams can be extended to a Dynkin diagram by an arm of length $k-1$
\cite{bfm}. Actually this has to be taken with a grain of salt: the generic codimension two 
linear section of the Grassmannian $G(2,2m+1)$ is also locally rigid, while according
to the previous rule, it would come from the diagram $E_{2m+2}$, which is not Dynkin 
if $m>3$! But these are the only exceptions.

\section{Some remarkable varieties with trivial canonical bundle} 

\subsection{More crepant resolutions}
Resolutions of singularities of orbit closures by total spaces of homogeneous 
vector bundles on generalized 
flag varieties are instances of  {\it Kempf collapsings}; we have already met three of them in
Example 4. 

The general construction is extremely 
simple. Let $E$ be a homogeneous
bundle on $G/P$, whose dual bundle is globally generated, and let 
$V^\vee:=H^0(G/P,E^\vee)$. Then $Tot(E)$ embeds inside $G/P\times V$, and the projection
to $V$ is the Kempf collapsing of $E$:
$$\xymatrix{
Tot(E)\ar@{->}[d]\ar@{->}[rr] & & Y\subset V\\
X & & 
}$$
These collapsings have nice properties, in
particular when the bundle $E$ is irreducible \cite{kempf, Weyman2003}.
They are always proper, so the image $Y$ is a closed subset of $V$, usually
singular, and the Kempf collapsing often provides a nice resolutions of singularities. 

Sometimes, these resolutions are even  {\it crepant}. One can show that this happens 
exactly when $\det(E)=\omega_X$ \cite{odl1}. Here are a few examples:

\medskip\noindent {\it Determinantal loci}. Let $V,W$ be vector spaces of dimensions $v,w$, and consider inside $Hom(V,W)$ the subvariety $D_k$ of morphisms of rank at most $k$. This means that the image is contained in (generically equal to) a $k$-dimensional subspace $U$ of $W$.
So $D_k$ is the image of a Kempf collapsing 
$$\xymatrix{
Tot(Hom(V,U_k))\ar@{->}[d]\ar@{->}[rr]& & D_k\subset Hom(V,W) \\
G(k,W) & & 
}$$
where $U_k$ is the rank $k$ tautological bundle on the Grassmannian. 
This construction is essential for constructing minimal resolutions of determinantal
loci \cite{Weyman2003}. This Kempf collapsing is crepant if and only if $V$ and $W$ have the same dimension, that is, for square matrices, independently of $k$. 

\medskip\noindent {\it Cubic polynomials}. Inside the spaces $C_k$ of cubic polynomials in 
$k$ variables, consider the set $C^2_k$ of those that can be written as polynomials in 
only two variables. This is the image of the Kempf collapsing 
$$\xymatrix{
Tot(S^3U_2)\ar@{->}[d]\ar@{->}[rr]^{p_1} & & C^2_k\subset C_k \\
G(2,k) & & 
}$$
Since $\det (S^3U_2)=\det (U_2)^6$, the morphism $p_1$ is crepant if and only if $k=6$.

\medskip\noindent {\it Skew-symmetric three-forms}. Inside the spaces $F_k$ of skew-symmetric three-forms 
in $k$ variables, consider the set $F^6_k$ of those that can be written as three-forms in 
only six variables. This is the image of the Kempf collapsing 
$$\xymatrix{
Tot(\wedge^3U_6)\ar@{->}[d]\ar@{->}[rr]^{p_2} & & F^6_k\subset F_k \\
G(6,k) & &
}$$
Since $\det (\wedge^3U_6)=\det (U_6)^{10}$, the morphism $p_2$ is crepant if and only if $k=10$.

\medskip\noindent {\it Pencils of quadrics}. Inside the spaces $P_k$ of pencils 
of quadrics in $k$ variables, consider the set $P_k^\ell$ of those that can be written in terms of $\ell$ variables only. This is the image of the Kempf collapsing 
$$\xymatrix{
Tot(\CC^2\otimes Sym^2U_\ell)\ar@{->}[d]\ar@{->}[rr]^{p_3} & & P^\ell_k\subset P_k \\
G(\ell,k) & & 
}$$
Since $\det (\CC^2\otimes Sym^2U_\ell)=\det (U_\ell)^{2\ell+2}$, the morphism $p_3$ is crepant if and only if $k=2\ell+2$. 

\subsection{Beauville-Donagi type constructions} Crepancy has the following nice 
consequence. 
Suppose that the crepancy condition is fulfilled for the vector bundle $E$ on $G/P$,
that is, $\omega_{G/P}=\det (E)$. Then by adjunction, a general
section $s$ of $E^\vee$ will vanish on a smooth subvariety $Z(s)\subset G/P$ with trivial
canonical bundle. 

Let us revisit the previous three examples from this perspective. 

\medskip\noindent {\it Cubic polynomials}. A global section $s$ of $S^3U_2^\vee$ on 
$G(2,6)$ is a cubic polynomial in six variables, and defines a cubic fourfold 
$X\subset \PP^5$. The variety $Z(s)\subset G(2,6)$ is the variety of lines $F(X)$ 
on this cubic fourfold. This is another fourfold, with trivial canonical bundle,
and Beauville and Donagi proved that it is hyperK\"ahler (\cite{bd}, more on this below).

\medskip\noindent {\it Skew-symmetric three-forms}. A global section $s$ of $\wedge^3U_6^\vee$ on 
$G(6,10)$ is an alternating three-form in ten variables, which also defines a hypersurface $X$ 
in $G(3,10)$. Then $Z(s)\subset G(6,10)$ can be interpreted as the parameter space for the 
copies of $G(3,6)$ contained in $X$. This is again a fourfold, with trivial canonical bundle,
and Debarre and Voisin proved that it is hyperK\"ahler \cite{dv}. 
Note the strong analogies with the 
previous case. 

\smallskip\noindent {\it Remark}. These two examples might give the impression that it 
should be easy to construct hyperK\"ahler manifolds as zero loci of sections of homogeneous
bundles on (generalized) Grassmannians. This impression is completely false. For instance,
Benedetti proved that if a hyperK\"ahler fourfold can be described as the zero locus of a 
general section of a semisimple homogeneous bundle on a Grassmannian, then it has to be one 
of the two previous examples \cite{vlad}. 

\medskip\noindent {\it Pencils of quadrics}. Consider a pencil $P\simeq\PP^1$ 
of quadrics in $k$ variables. In general it contains exactly $k$ singular quadrics, which
have corank one. For $k=2\ell+2$ even, the family of $\ell$-dimensional projective spaces that 
are contained in some quadric of the pencil defines a hyperelliptic curve $C\stackrel{\eta}{\lra} P$ of genus $\ell$. Moreover, the family of $(\ell-1)$-dimensional projective spaces that 
are contained in the base locus of the pencil is a $\ell$-dimensional smooth manifold $A$ 
with trivial canonical bundle. Reid proved that $A$ is an abelian variety, isomorphic
with the Jacobian variety of $C$ \cite[Theorem 4.8]{reid}. 

\subsection{The Springer resolution and nilpotent orbits}
A fundamental homogeneous vector bundle on $G/P$ is the cotangent bundle $\Omega_{G/P}$. Its  
dual, the tangent bundle, is generated by global sections. Moreover, it is a general
fact that the space of global sections of the tangent bundle is the Lie algebra
of the automorphism group. Beware that the automorphism group of $G/P$ can very-well 
be bigger than $G$ (we have seen that for $\QQ^5=G_2/P_1=B_3/P_1$), but except for 
a few well understood exceptions, we have $H^0(G/P,T_{G/P})=\fg$. 

In any case, $G/P$ being $G$-homogeneous, there is always an injective map from $\fg$
to $H^0(G/P,T_{G/P})$, whose image is a linear system of sections that generate the 
tangent bundle at every point. In particular the Kempf collapsing of $\Omega_{G/P}$
to $\fg$ is well-defined.
Note moreover that the fiber of $T_{G/P}$ over the base point is $\fg/\fp$, 
whose dual identifies with $\fp^\perp\subset\fg^\vee\simeq\fg$ (where the duality is given 
by the Cartan-Killing form). One can check
that $\fp^\perp$ is the nilpotent radical of $\fp$, in particular it is made of 
nilpotent elements of $\fg$. 

One denotes by $\cN$ the {\it nilpotent cone} in $\fg$. 
A fundamental result in Lie theory is that the adjoint action of $G$ on $\cN$ has finitely
many orbits. This implies that for each parabolic $P$, the image of the Kempf 
collapsing of the cotangent bundle of $G/P$ is the closure of a uniquely defined 
nilpotent orbit $\cO_P$:
$$\xymatrix{
Tot(\Omega_{G/P})\ar@{->}[d]\ar@{->}[rr]^{\pi_P} & & \bar\cO_P\subset \cN\subset \fg \\
G/P & & 
}$$
Note that the crepancy condition is automatically fulfilled. 
In fact both sides of the collapsing have natural symplectic structures. These are preserved by
$\pi_P$, but beware that this collapsing is not neccessarily birational, although it 
is in most cases. The most important one is the case where $P=B$, which yields the 
so-called {\it Springer resolution} \cite{cg}
$$\xymatrix{
Tot(\Omega_{G/B})\ar@{->}[d]\ar@{->}[rr]^{\pi_B} & & \cN\subset \fg \\
G/B & & 
}$$
This is a fundamental example of a symplectic resolution of singularities. 

\subsection{Stratified Mukai flops}
In some situations the same orbit closure has several distinct symplectic resolutions.
This happens in particular when the Dynkin diagram of $G$ has a non trivial symmetry, 
that is in type $A, D$ or $E_6$. In type $A$, the cotangent bundles of the Grassmannians
$G(k,n)$ and $G(n-k,n)$ resolve the same nilpotent orbit closure, so there is a birational
morphism between them which is called a {\it stratified Mukai flop} (of type $A$). 
$$\xymatrix{ & F(k,n-k,n)\ar@{->}[ld]\ar@{->}[rd] & \\
G(k,n) & Tot(E_0)\ar@{->}[ld]\ar@{->}[u]\ar@{->}[rd] & G(n-k,n) \\
 Tot(\Omega_{G(k,n)})\ar@{->}[rd]\ar@{->}[u] && 
 Tot(\Omega_{G(n-k,n)})\ar@{->}[ld]\ar@{->}[u] \\
 & \bar\cO_{k,n} &
}$$
where $\bar\cO_{k,n}\subset\fsl_n$ is the set of nilpotent matrices of square zero and 
rank at most $k$. In particular $k=1$ gives a {\it Mukai flop}. 
The other cases yield the so-called stratified Mukai flops of types D and E.
In particular the cotangent spaces of the two spinor varieties of type $D_n$ resolve the 
same nilpotent orbit \cite{nam03, fu}.

\subsection{Projective duality}
Consider an embedded projective variety $X\subset\PP (V)$ and its projective dual $X^*\subset \PP (V^\vee)$. 
By definition 
$X^*$ parametrizes the tangent hyperplanes to $X$. If the latter is smooth, this implies 
that the affine cone over $X^*$ is the image of the Kempf collapsing 
$$\xymatrix{
Tot((\hat TX)^\perp )\ar@{->}[d]\ar@{->}[rr]^{\pi} & & \hat X^*\subset V^\vee\\
X & & 
}$$
where we denoted by $\hat TX$ the affine tangent bundle of $X$, which is a sub-bundle of 
the trivial bundle with fiber $V$.
If $\pi$ is generically finite, then $X^*\subset \PP (V^\vee)$ is a hypersurface: this is 
the general expectation. Moreover projective duality is an involution, and this implies that
the general fibers of $\pi$ are projective spaces. In particular, if $X^*$ is indeed a 
hypersurface, then the projectivization of $\pi$ is a resolution of singularities. 

\smallskip Now suppose that $V$ and its dual are $G$-modules with finitely many orbits. 
Since projective duality is compatible with the $G$-action, it must define a bijection 
between the $G$-orbit closures in $\PP(V)$ and those in $\PP(V^\vee)$ (but not necessarily 
compatible with inclusions). 

\subsection{Cubic-K3 and Pfaffian-Grassmannian dualities}
To be even more specific, consider $G=GL_n$ acting on the space 
$V=\wedge^2(\CC^n)^\vee$ 
of alternating bilinear forms. The orbits are defined by the rank, which can be any even integer
between $0$ and $n$. We denote by $Pf_r$ the orbit closure consisting of forms of rank at most $r$.
These orbits can also be seen inside $V^\vee$. One has 
$$(Pf_r)^*=Pf_{n-r-\epsilon},$$
with $\epsilon=0$ if $n$ is even, $\epsilon=1$ if $n$ is odd. For $n$ even, $Pf_{n-2}$ is the Pfaffian hypersurface, of degree $n/2$. But if $n$ is odd, there is no invariant hypersurface, 
and the complement $Pf_{n-3}$ of the open orbit has codimension three. 

\subsubsection*{$n=5$} The $PGL_5$-orbits in $\PP (\wedge^2\CC^5)$ are 
$G(2,5)$ and its complement. In particular $G(2,5)$ is projectively self-dual. Note that its 
index is $5$, so its intersection 
$$X_g=G(2,5)\cap gG(2,5)$$ with a translate by  a general $g\in PGL(\wedge^2\CC^5)$ 
is a smooth Calabi-Yau threefold. This is also true for the intersection
of their projective duals, 
$$Y_g=G(2,5)^*\cap (gG(2,5))^*.$$
These two Calabi-Yau's are obviously deformation equivalent, but \cite{bcp, or}: 

\begin{theorem}
For $g\in PGL(\wedge^2\CC^5)$ general: 
\begin{enumerate}
\item $X_g$ and $Y_g$ are derived equivalent.
\item $X_g$ and $Y_g$ are not birationally equivalent.
\item In the Grothendieck ring of varieties, $([X_g]-[Y_g])\times \LL^4 = 0.$ 
\end{enumerate}
\end{theorem}

Note moreover that Example 2 is a limit case, obtained when the two Grassmannians collapse one to the other along some normal direction (indeed the 
normal bundle to $G(2,5)\subset \PP^9$ is isomorphic with $Q^\vee(2)$). Interestingly, 
for this case one obtains in the Grothendieck ring the stronger relation 
$([Z_1(s)]-[Z_2(s)])\times \LL^2 = 0.$

Similar phenomena can be observed if one replaces the Grassmannian $G(2,5)$ by the spinor
variety $\SS_{10}$ \cite{manspin}.

\subsubsection*{$n=6$} The proper $PGL_6$-orbit closures  in $\PP (\wedge^2\CC^6)$ are the Grassmannian $G(2,6)$ and the Pfaffian hypersurface. The projective dual of $G(2,6)$ is the 
Pfaffian hypersurface $Pf_4$ inside $\PP (\wedge^2\CC^6)^\vee\simeq \PP (\wedge^4\CC^6)$, 
which is singular exactly along $G(4,6)$. Since the codimension of the latter is $6$, the 
intersection of $Pf_4$ with a general five dimensional projective space $L\subset \PP (\wedge^4\CC^6)$
is a smooth cubic fourfold $X_L$. Moreover $L^\perp\subset \PP (\wedge^2\CC^6)$ has codimension 
$6$, and it follows that its intersection with $G(2,6)$ is a smooth K3 surface $S_L$. 
Beauville and Donagi \cite{bd} showed that:

\begin{proposition} The Fano variety of lines $F(X_L)\simeq Hilb^2(S_L)$. 
\end{proposition}

In particular $F(X_L)$ is hyperK\"ahler, and since this property is preserved by deformation, 
the variety of lines of any cubic fourfold is also hyperK\"ahler, as soon as it is smooth.

\proof Take two points in $S_L$, defining two planes $P_1$ and $P_2$ in $\CC^6$. In general 
these planes are transverse and we can consider $Q=P_1\oplus P_2$. The alternating form that
vanish on $Q$ (which implies they are degenerate) span a $\PP^5$ in $\PP (\wedge^2\CC^6)^\vee$, 
which is automatically orthogonal to $P_1$ and $P_2$. So the intersection with $L$ is defined
by only four conditions (rather than $6$), and the intersection is in general a line in $X_L$.
This defines a rational map $a: Hilb^2(S_L)\lra F(X_L)$, which is in fact regular. 

Conversely, take a line $\ell$ in $X_L$. In particular this is a line parametrizing 
alternating forms
of rank four, and if the rank drops to two $X_L$ must be singular, which we exclude. By 
\cite[Proposition 2]{mm}, we have:

\begin{lemma} 
Up to the action $GL_6$, a line of alternating forms of constant rank four in six variables can
be put in one of the two possible normal forms
$$\langle e_1\wedge e_3+e_2\wedge e_4, e_1\wedge e_5+e_2\wedge e_6\rangle, \quad \mathit{or}\quad
\langle e_1\wedge e_3+e_2\wedge e_4, e_1\wedge e_4+e_2\wedge e_6\rangle .$$
\end{lemma}

In particular, there always exists a unique four plane $Q$ (the orthogonal to $e_1$ and $e_2$) 
on which all the forms in the line $\ell$ vanish. Note that $\PP(\wedge^2Q)^\perp\cap L$ is 
a linear space in $X_L$ that contains $\ell$, so it must coincide with $\ell$ unless $X_L$
contains a plane. By dimension count, this implies that $\PP(\wedge^2Q)\cap L^\perp$ is also
a line, which has to cut the quadric $G(2,Q)\subset \PP(\wedge^2Q)$ along two points (more
precisely along a length two subscheme). These two points belong to $S_L$, so this 
defines a map $b: F(X_L)\lra Hilb^2(S_L)$, which is inverse to $a$. \qed

\subsubsection*{$n=7$} The proper $PGL_7$-orbit closures in $\PP (\wedge^2\CC^7)$ are the Grassmannian $G(2,7)$ and the Pfaffian variety of degenerate forms, of codimension three.
The projective dual of $G(2,7)$ is the 
Pfaffian variety $Pf_4$ inside $\PP (\wedge^2\CC^7)^\vee\simeq \PP (\wedge^5\CC^7)$, 
which is singular exactly along $G(5,7)$. The codimension of the latter is $10$, so the 
intersection of $Pf_4$ with a general six dimensional projective space $L\subset \PP (\wedge^5\CC^7)$ is a smooth threefold $X_L$. Moreover $L^\perp\subset \PP (\wedge^2\CC^7)$ has codimension $7$, so its intersection with $G(2,7)$ is a smooth Calabi-Yau threefold $Y_L$. 

\begin{theorem}\label{bc} Suppose that $X_L$ and $Y_L$ are smooth. Then: 
\begin{enumerate}
\item $X_L$ and $Y_L$ are derived equivalent Calabi-Yau threefolds.
\item $X_L$ and $Y_L$ are not birationally equivalent. 
\item In the Grothendieck ring of varieties, $([X_L]-[Y_L])\times \LL^6=0$. 
\end{enumerate}
\end{theorem}

This is the famous {\it Pfaffian-Grassmannian equivalence} \cite{bc}. The identity in the Grothendieck ring is due to Martin. 

\proof The fact that $X_L$ has trivial canonical bundle can be proved by observing 
that the minimal resolution of the Pfaffian locus $Pf_4\subset\PP := \PP (\wedge^5\CC^7)$ is 
the Buchsbaum-Eisenbud complex
$$0\lra \cO_\PP(-7)\stackrel{\wedge^3\omega}{\lra} \cO_\PP(-4)\stackrel{\omega}{\lra}
 \cO_\PP(-3)\stackrel{\wedge^3\omega}{\lra} \cO_\PP\lra \cO_{Pf_4}\lra 0.$$
This implies that the relative canonical bundle of $Pf_4\subset \PP$ is the restriction of 
$\cO(-7)$, and this remains true by taking a linear section. So the canonical bundle 
of $X_L$ must be trivial. Using the resolution of $\cO_{X_L}$ obtained by restricting 
to $L$ the previous complex, one computes that $h^1(\cO_{X_L})=0$, and this is enough
to ensure that $X_L$ is indeed Calabi-Yau. 

The same argument as in the proof of Theorem 2.3 implies that $X_L$ and $Y_L$ are not 
birationally equivalent. Proving that they are derived equivalent requires more 
sophisticated techniques. \qed

\smallskip \noindent {\it Remark}. The Calabi-Yau threefolds in the $G_2$-Grassmannians 
that appear in the proof of Theorem \ref{zerodiv} are degenerations of the Pfaffian-Grassmannians
Calabi-Yau threefolds. 

\medskip
All these examples are in fact instances of Kuznetsov's Homological Projective Duality, 
which provides a way, under favourable circumstances, to compare the derived category of 
a variety $X\subset\PP (V)$, and of its linear sections, with the derived category of the
dual $X^*\subset \PP (V^\vee)$, and of its dual linear sections (see \cite{thomas} for an introduction). 

\smallskip Another example with a similar flavor is that of the spinor variety $\SS_{10}\subset\PP(\Delta_{16})$. As we already mentionned, 
this variety is a very similar to $G(2,5)$. 
In particular it is projectively self-dual (non canonically), and its complement 
is an orbit of $Spin_{10}$. Moreover it is a Mukai variety, so a smooth codimension 
eight linear section $X_M=\SS_{10}\cap M$ is a K3 surface of degree $12$. Moreover $M^\perp$
also has codimension eight and we get another K3 surface $Y_M=\SS_{10}^*\cap M^\perp$
in the dual projective space. These two K3 surfaces are Fourier-Mukai partners (they
are derived equivalent), they are not isomorphic, and $([X_M]-[Y_M])\times \LL = 0$ 
in the Grothendieck ring of varieties \cite{imou16}.

\section{Orbital degeneracy loci} 

Instead of zero loci of sections of vector bundles, it is very natural to consider morphisms
between vector bundles and their degeneraly loci, where the morphism drops rank. More generally,
one may consider {\it orbital degeneracy loci} and try to construct interesting varieties 
from those. 

\subsection{Usual degeneracy loci} Given a morphism $\phi : E\lra F$ 
between vector bundles on a variety $X$, the degeneracy loci are the closed subvarieties
where the rank drops:
$$D_k(\phi):=\{x\in X, \; rk (\phi_x : E_x\lra F_x)\le k\}.$$
Writing locally $\phi$ as a matrix of regular functions and taking $(k+1)$-minors allow to give
to $D_k(\phi)$ a canonical scheme structure.
If $\phi$ is sufficiently regular, $D_k(\phi)$ has codimension $(e-k)(f-k)$ (where $e$ and
$f$ are the ranks of the bundles $E$ and $F$), and its singular locus is exactly $D_{k+1}(\phi)$.
If moreover the latter is empty, the kernel and cokernel of $\phi$ give well-defined vector 
bundles $K$ and $C$ on $D_k(\phi)$, with a long exact sequence 
$$0\lra K\lra E_{D_k(\phi)}\lra F_{D_k(\phi)}\lra C\lra 0.$$
A local study allows to check that the normal bundle $N_{D_k(\phi)/X}\simeq Hom(K,C)$. Using the 
previous exact sequence, we can deduce in the square format $e=f$ that the relative canonical
bundle is 
$$\omega_{D_k(\phi)/X}\simeq (\det E^\vee)^k\otimes (\det F)^k.$$
This gives a nice way to construct Calabi-Yau manifolds of small dimension: start from a Fano variety $X$ such that $\omega_X=L^{-k}$ for some ample line bundle $L$; find vector bundles 
$E, F$ of the same rank $e$ such that $Hom(E,F)$ is globally generated and $(\det E^\vee)\otimes (\det F)=L$. If the dimension of $X$ is $n<(e-k+1)^2$, a general element $\phi\in Hom(E,F)$ 
defines a smooth manifold $D_k(\phi)$ with trivial canonical bundle (most often a Calabi-Yau), of 
dimension $n-(e-k)^2$. Similar remarks apply to the case of skew-symmetric morphisms 
$\phi : F^\vee\lra F$ and their associated Pfaffian loci.

\smallskip\noindent {\it Example 5}. Consider a three-form $\omega\in\wedge^3(\CC^n)^\vee$, 
and recall that it defines a global section of the bundle $\wedge^2Q^\vee (1)$ on $\PP^{n-1}$. 
Locally this is just a family of two-forms on the quotient bundle, and we get a stratification 
of $\PP^{n-1}$ by Pfaffian loci.

\subsection{General theory} Determinantal and Pfaffian loci are special cases of a more 
general construction, that of orbital degeneracy loci. Suppose for example that 
$E$ is a vector bundle on some variety $X$, with rank $e$. Consider a closed subvariety
$Y$ of $\wedge^m\CC^e$, stable under the action of $GL_e$. Let $s$ be a section of $\wedge^mE$. 
Then the locus of points $x\in X$ for which $s(x)\in\wedge^mE_x\simeq \wedge^m\CC^e$ belongs 
to $Y$ is well-defined, since the latter identification, although non canonical, only varies
by an element of $GL_e$. This locus is precisely the orbital degeneracy locus $D_Y(s)$. 

Basic properties of these loci are established in \cite{odl1}. For instance, for a general
section $s$ one has
$$\mathrm{codim}(D_Y(s),X)=\mathrm{codim}(Y, \wedge^m\CC^e),$$
$$\mathrm{Sing}(D_Y(s))=D_{Sing(Y)}(s).$$

\subsection{Applications: Coble cubics and generalized Kummers}
For a nice example of these constructions, start with  a general three-form $\omega\in\wedge^3(\CC^9)^\vee$. The quotient bundle $Q$ on $\PP^8$ has rank $8$, 
and a skew-symmetric form in eight variables has rank six in codimension $1$, four in 
codimension $6$, two in codimension $15$. So the Pfaffian stratification reduces
to $\PP^8\supset C\supset A$, where:
\begin{enumerate}
\item The hypersurface $C$ is given by the Pfaffian of the induced section of $\wedge^2Q^\vee(1)$, 
which is a degree four polynomial, hence a section of the line bundle $\wedge^8Q^\vee (4)=\cO(3)\subset Sym^4(\wedge^2Q^\vee(1))$; so $C$ is in fact a cubic hypersurface. 
\item The singular locus of $C$ is the smooth surface $A$, whose canonical bundle is trivial. 
This is in fact an {\it abelian} surface, and $C$ is the {\it Coble cubic} of $A$, the unique cubic hypersurface which is singular exactly along $A$ \cite{gsw}.
\end{enumerate}
On the dual projective space, denote by $H$ the tautological hyperplane bundle, of rank $8$.
Then $\omega$ defines a global section of $\wedge^3H^\vee$, and therefore, orbital 
degeneracy loci $D_Y(\omega)$ for each orbit closure $Y\subset\wedge^3\CC^8$.
The codimension four orbit closure $Y_4$ has particular interest \cite{coble}. It can be characterized as the image of the birational Kempf collapsing 
$$\xymatrix{
Tot(E )\ar@{->}[d]\ar@{->}[rr] & & Y_4\subset \wedge^3\CC^8 \\
F(2,5,8) & & 
}$$
where $E$ is the rank $31$ vector bundle with fiber $\wedge^3U_5+U_2\wedge U_5\wedge \CC^8$ over the partial flag $U_2\subset U_5\subset \CC^8$.

\begin{theorem}
The Kempf resolution of the orbital degeneracy locus $D_{Y_4}(\omega)$ is a hyperK\"ahler fourfold.

More precisely, this   hyperK\"ahler fourfold is isomorphic with the generalized Kummer fourfold of the 
abelian surface $A$. 
\end{theorem}

Generalized Kummers of abelian surfaces were first constructed
by Beauville: start with an abelian surface $A$ and the Hilbert scheme $Hilb^n(A)$ of 
$n$-points on $A$; then the sum map $A^n\ra A$ descends to $Hilb^n(A)\ra A$, and 
the $n$-th generalized Kummer variety of $A$ is any fiber. The resolution in Theorem 4.1 
is exactly the restriction to the Kummer variety of the Hilbert-Chow map from $Hilb^3(A)$
to $Sym^3(A)$. 

Generalized Kummer varieties are 
hyperK\"ahler manifolds, like Hilbert schemes of points on K3 surfaces. But contrary to
the latter case (for which we have e.g. Fano varieties of lines on cubic fourfolds),
no locally complete projective deformation of generalized Kummers is known!

\end{document}